\documentclass[12pt]{article}
\usepackage{amssymb}
\usepackage{amsfonts}
\usepackage{amsmath}

\newtheorem{theorem}{Theorem}

\newtheorem{corollary}[theorem]{Corollary}

\newtheorem{definition}[theorem]{Definition}

\newtheorem{lemma}[theorem]{Lemma}

\newtheorem{proposition}[theorem]{Proposition}
\newtheorem{remark}[theorem]{Remark}

\newenvironment{proof}[1][Proof]{\noindent\textbf{#1.} }{\ \rule{0.5em}{0.5em}}
\input{tcilatex}

\begin{document}

\title{Dimension and hitting time in rapidly mixing systems}
\author{S. Galatolo \\
Dipartimento di Matematica Applicata, Via Buonarroti 1, Pisa\\
galatolo@dm.unipi.it, d80288@ing.unipi.it}
\maketitle

\begin{abstract}
We prove that if a system has superpolynomial (faster than any power law)
decay of correlations then the time $\tau _{r}(x,x_{0})$ needed for a
typical point $x$ to enter for the first time a ball $B(x_{0},r)$ centered
in $x_{0},$ with small radius \ $r$ scales as the local dimension at $x_{0},$
i.e.%
\begin{equation*}
\underset{r\rightarrow 0}{\lim }\frac{\log \tau _{r}(x,x_{0})}{-\log r}%
=d_{\mu }(x_{0}).
\end{equation*}

This result is obtained by proving a kind of dynamical Borel-Cantelli lemma
wich holds also in systems having polinomial decay of correlations.
\end{abstract}

\section{Introduction}

Let us consider an ergodic system $(X,T,\mu )$ and two typical points $%
x,x_{0}.$ The orbit of $x$ will come as near as we want to $x_{0}$ entering
(sooner or later) in each positive measure neighborhood of $x_{0}.$ If $%
S\subset X$ is such a neighborhood, there is $n\in \mathbb{N}$ such that $%
T^{n}(x)\in S.$ Hitting time (also called waiting time, or shrinking target)
problems consider the time 
\begin{equation*}
\tau _{S}(x)=\min \{n\geq 1:T^{n}(x)\in S\}
\end{equation*}%
needed for the orbit of $x$ to enter in $S$ for the first time. As the
measure $\mu (S)$ is smaller and smaller, $\tau _{S}(x)$ is bigger and
bigger. In systems whose behavior is chaotic enough it can be expected that
since (by ergodicity) the orbit of $x$ must visit $S$ with (asymptotic)
frequency $\mu (S),$ when $S$ is small, for most $x$%
\begin{equation}
\tau _{S}(x)\sim \frac{1}{\mu (S)}.\footnote{%
Here "$\sim $" stands for some kind of (more or less strict) equivalence in
the asymptotic behavior of the two quantities.}  \label{XXX}
\end{equation}%
In general systems it is not always like this, even for nice sets, such as
balls. There are ergodic systems where the waiting time is generically
larger than $\frac{1}{\mu (S)}.$ In some sense the behavior of such systems
converges "slowly" to the ergodic one.

A way to express this kind of problem in more precise terms is to consider a
sequence of balls $B(x_{0},r_{k})$ centered in $x_{0}$ and consider the
scaling behavior of the waiting time in such balls, considering as an
indicator for the waiting time $R(x,x_{0})=\underset{k\rightarrow \infty }{%
\lim }\frac{\log (\tau _{B(x_{o},r_{k})}(x))}{-\log (r_{k})}$ (hence $\tau
_{B(x_{o},r_{k})}(x)\sim $ $r_{k}^{-R(x,x_{0})}$ in the sense of scaling
behavior equivalence). On the other hand we can consider the scaling
behavior of the measure of the balls: $d_{\mu }(x_{0})=\underset{%
k\rightarrow \infty }{\lim }\frac{\log \mu (B(x_{0},r_{k}))}{\log (r_{k})}$
(as above this give $\mu (B(x_{0},r_{k}))\sim r_{k}^{d_{\mu }(x_{0})}$) this
is called the local dimension of $\mu $ at $x_{0}.$ It is not difficult to
prove that in general systems $R(x,x_{0})\geq d_{\mu }(x_{0})$ (for precise
statements see theorem \ref{GAN}) and there are systems (rotations by
irrational angles which are well approximable, see \cite{KS} e.g.) where $%
\underset{k\rightarrow \infty }{\lim \sup }\frac{\log (\mu (B(x_{0},r_{k})))%
}{\log (r_{k})}>d_{\mu }(x_{0})$. Conversely, in many systems with more or
less chaotic behavior we have the analog of equation \ref{XXX} that is%
\begin{equation}
R(x,x_{0})=d_{\mu }(x_{0})  \label{YYY}
\end{equation}%
for $\mu -$almost each $x.$ This relation was proved for example for Axiom A
systems (and in a weaker form for typical Interval Exchange Transformations,
see \cite{G2}), typical rotations on the circle (\cite{KS}), systems having
exponential distribution of return times (\cite{G}, theorem 6) or systems
having a strong mixing assumption (called uniform mixing) which is verified
in a class of complex dynamical systems (inner functions, see \cite{FMP}).

Hitting time problems are related to many other features of chaotic
dynamics: entropy (\cite{Sh}), quantitative recurrence and the distribution
of return times (see e.g. \cite{G}, \cite{LHV}), the so called dynamical
Borel-Cantelli property (see Definition \ref{sbc}, and references \cite{GK}
or \cite{FMP}), asymptotics of Birkhoff sums of functions having infinite
average (see \cite{G2}), orbit complexity (\cite{BGI}). Borel-Cantelli
properties have also relations with rate of mixing and speed of
approximations of points (see e.g. \cite{T}, \cite{K}, \cite{FMP}).

As is well known, in a mixing system we have $\mu (A\cap
T^{-n}(B))\rightarrow \mu (A)\mu (B)$ for each measurable sets $A,B.$ The
speed of convergence of the above limit can be arbitrarily slow (depending
on $T$ but also on the shape of the sets $A,B$). In many systems however the
speed of convergence can be estimated for sets having some regularity. Let
us remark that considering $1_{A}(x)=\left\{ 
\begin{array}{c}
1~if~x\in A \\ 
0~if~x\notin A%
\end{array}%
\right. $ the mixing condition becomes $\int 1_{B}\circ T^{n}1_{A}d\mu
\rightarrow \int 1_{A}d\mu \int 1_{B}d\mu .$ If we consider Hoelder
observables $\phi ,\psi :X\rightarrow \mathbb{R}$, in many systems (most of
them having some form of hyperbolicity) it can be proved that 
\begin{equation*}
|\int \phi \circ T^{n}\psi d\mu -\int \phi d\mu \int \psi d\mu |\leq
\left\vert \left\vert \phi \right\vert \right\vert \left\vert \left\vert
\psi \right\vert \right\vert \Phi (n)
\end{equation*}%
where $||\ast ||$ is some Hoelder norm, and $\Phi (n)\rightarrow 0$ is a
function whose decay rate is estimated. For example, considering axiom A
systems with equilibrium states we have a similar result and $\Phi (n)$
decays exponentially fast. Similar results have been obtained for a large
class of more or less hyperbolic systems. The set of references for this
kind of results is huge, we cite the books \cite{B}, \cite{V} and \cite{L},
a survey of recent results on nonuniformly expanding systems. Recently fast
decay of correlation has also been proved in a large class of systems having
strange attractors with dimension close to 1 (Rank 1 systems, see \cite{WY}, 
\cite{WY2}). This class, contains for example the Henon map for some
interesting set of parameters.

In this paper we prove (theorem \ref{maine}) that if the system has a decay
of correlations (definition \ref{sup}) faster than any power law, then for 
\emph{each }$x_{0}\in X$ such that $d_{\mu }(x_{0})$ exists, eq. \ref{YYY}
holds for almost each $x$. Hence for this kind of systems the typical
hitting time scaling rate is given by the local dimension at $x_{0}$.

As a corollary of the proof (corollary \ref{XXX2}) we obtain that in systems
with fast enough decay of correlations (including some case of polynomial
decay) a large class of decreasing sequences of sets have the strong
Borel-Cantelli property (see definition \ref{sbc}). About this, we remark
that in \cite{T} and \cite{FMP} other Borel-Cantelli results are obtained
supposing different assumptions on the system's mixing behavior (in most
cases it is supposed to have some kind of uniform bound on the speed of
mixing).

We end by remarking that since eq. \ref{YYY} is proved for each $x_{0}$
(having local dimension) the hitting time can be used to numerically
estimate local dimension at \emph{nontypical} points. For an example of
numerical use of hitting times for this kind of question see \cite{CG}.

\section{Dimension and waiting time}

In the following we will consider a discrete time dynamical system $(X,T)$
where $X$ is a separable metric space equipped with a Borel finite measure $%
\mu $ and $T:X\rightarrow X$ is a measurable map.

Let us consider the first entrance time of the orbit of $x$ in the ball $%
B(x_{0},r)$ with center $x$ and radius $r$ 
\begin{equation*}
\tau _{r}(x,x_{0})=\min (\{n\in \mathbf{N},n>0,T^{n}(x)\in B(x_{0},r)\})\,.
\end{equation*}%
By considering the power law behavior of $\tau _{r}(x,x_{0})$ as $%
r\rightarrow 0$ let us define the hitting time indicators as%
\begin{equation*}
\overline{R}(x,x_{0})=\mathrel{\mathop{limsup}\limits_{r\rightarrow 0}}\frac{%
\log (\tau _{r}(x,x_{0}))}{-\log (r)},\underline{R}(x,x_{0})=%
\mathrel{\mathop{liminf} \limits_{r\rightarrow 0}}\frac{\log (\tau
_{r}(x,x_{0}))}{-\log (r)}.
\end{equation*}

If for some $r,$ $\tau _{r}(x,x_{0})$ is not defined then $\overline{R}%
(x,x_{0})$ and $\underline{R}(x,x_{0})$ are set to be equal to infinity. The
indicators $\overline{R}(x)$ and $\underline{R}(x)$ of quantitative
recurrence defined in \cite{BS} are obtained as a special case, $\overline{R}%
(x)=\overline{R}(x,x)$, $\underline{R}(x)=\underline{R}(x,x)$.

We recall some basic properties of $R(x,x_{0})$ which follow from the
definition:

\begin{proposition}
\label{inizz}$R(x,x_{0})$ satisfies the following properties

\begin{itemize}
\item $\overline{R}(x,x_{0})=\overline{R}(T(x),x_{0})$, $\underline{R}%
(x,x_{0})=\underline{R}(T(x),x_{0})$.

\item If $T$ is $\alpha -Hoelder$, then $\overline{R}(x,x_{0})\geq \alpha 
\overline{R}(x,T(x_{0}))$, $\underline{R}(x,x_{0})\geq \alpha \underline{R}%
(x,T(x_{0}))$.
\end{itemize}
\end{proposition}

Before proceeding and stating connections between hitting time and local
dimension, let us recall more precisely some results about dimension. If $X$
is a metric space and $\mu $ is a measure on $X$ the local dimension of $\mu 
$ at $x$ is defined as $d_{\mu }(x)=\mathrel{\mathop{lim}\limits_{r%
\rightarrow 0}}\frac{log(\mu (B(x,r)))}{log(r)}$ (when the limit exists).
Conversely, the upper local dimension at $x\in X$ is defined as $\overline{d}%
_{\mu }(x)=\mathrel{\mathop{limsup}\limits_{r\rightarrow 0}}\frac{log(\mu
(B(x,r)))}{log(r)}$ and the lower local dimension $\underline{d}_{\mu }(x)$
is defined in an analogous way by replacing $limsup$ with $liminf$. If $%
\overline{d}_{\mu }(x)=\underline{d}_{\mu }(x)=d$ almost everywhere the
system is called exact dimensional. In this case many notions of dimension
of a measure will coincide. In particular $d$ is equal to the dimension of
the measure: $d=\inf \{\dim _{H}Z:\mu (Z)=1\}.$ This happens in a large
class of systems. For example in systems having nonzero Lyapunov exponents
almost everywhere (see for example the book \cite{P}).

In general systems the quantitative recurrence indicator gives only a \emph{%
lower} bound on the dimension. The hitting time indicator instead gives an 
\emph{upper} bound to the local dimension of the measure at the point $y$.
This is summarized in the following

\begin{theorem}
\label{GAN}(\cite{BGI},\cite{G},\cite{BS}, \cite{Bo}) If $(X,T,\mu )$ is a
dynamical system over a separable metric space, with an invariant measure $%
\mu .$ For each $x_{0}$ 
\begin{equation}
\underline{R}(x,x_{0})\geq \underline{d}_{\mu }(x_{0})\ ,\ \overline{R}%
(x,x_{0})\geq \overline{d}_{\mu }(x_{0})  \label{thm4}
\end{equation}%
holds for $\mu $ almost each $x$. Moreover, if $X$ is a closed subset of ${%
\mathbb{R}
}^{n}$, then for almost each $x\in X$ 
\begin{equation}
\overline{R}(x,x)\leq \overline{d}_{\mu }(x)\ ,\ \underline{R}(x,x)\leq 
\underline{d}_{\mu }(x)\,.  \label{Bsa}
\end{equation}
\end{theorem}

We remark that eq. \ref{thm4} implies that if $\alpha <\underline{d}_{\mu
}(x_{0})$ then 
\begin{equation*}
\underset{n\rightarrow \infty }{\lim \inf }n^{\frac{1}{\alpha }%
}d(T^{n}(x),x_{0})=\infty .
\end{equation*}%
In systems with superpolynomial decay of correlations the above inequalities
become equalities. For eq. \ref{Bsa} this is proved in \cite{S}. We are
going to consider eq. \ref{thm4}. Let us recall more precisely what is
superpolynomial decay of correlations.

\begin{definition}
\label{sup}Let $\phi ,$ $\psi :X\rightarrow \mathbb{R}$ be Lipschitz
observables on $X$. A system $(X,T,\mu )$ is said to have superpolynomial
decay of correlations if 
\begin{equation*}
|\int \phi \circ T^{n}\psi d\mu -\int \phi d\mu \int \psi d\mu |\leq
\left\vert \left\vert \phi \right\vert \right\vert \left\vert \left\vert
\psi \right\vert \right\vert \Phi (n)
\end{equation*}%
with $\Phi $ having superpolynolmial decay, i.e. $\lim n^{\alpha }\Phi
(n)=0, $ $\forall \alpha >0.$ Here $||~||$ is the Lipschitz norm\footnote{%
The case where the Holder norms are considered follows from the case of
Lipschitz norm.}.
\end{definition}

With this definition we can state the main result of the paper:

\begin{theorem}
\label{maine} If $(X,T,\mu )$ has superpolynomial decay of correlations and $%
d_{\mu }(x_{0})$ exists then%
\begin{equation}
\overline{R}(x,x_{0})=\underline{R}(x,x_{0})=d_{\mu }(x_{0})  \label{maineq}
\end{equation}%
for $\mu $-almost each $x$.
\end{theorem}

We remark that eq. \ref{maineq} easily implies that if $\alpha >\underline{d}%
_{\mu }(x_{0})$ then 
\begin{equation*}
\underset{n\rightarrow \infty }{\lim \inf }n^{\frac{1}{\alpha }%
}d(T^{n}(x),x_{0})=0.
\end{equation*}%
Before proving theorem \ref{maine} we need some lemmas. A sequence of sets $%
S_{n}\subset X$ is said to be strongly Borel-Cantelli if in some sense the
preimages $T^{-n}S_{n}$ cover the space uniformly:

\begin{definition}
\label{sbc}Let $1_{S}$ be the indicator function of the set $S.$ The
sequence of subsets $S_{n}\subset X$ is said to be a Strongly Borel-Cantelli
sequence (SBC) if $\sum_{n}\mu (S_{n})=\infty $ and for $\mu -$a.e. $x\in X$
we have as $N\rightarrow \infty $ 
\begin{equation*}
\frac{\sum_{n=1}^{N}1_{T^{-n}S_{n}}(x)}{\sum_{n=1}^{N}\mu (S_{n})}%
\rightarrow 1,~\mu -a.e.
\end{equation*}
\end{definition}

We remark that posing $Z_{k}(x)=\sum_{0}^{k}1_{T^{-i}S_{i}}(x)$ the above
condition is equivalent to 
\begin{equation*}
\frac{Z_{k}}{E(Z_{k})}\rightarrow 1,~\mu -a.e.
\end{equation*}%
The following technical lemma estimates the speed of mixing of balls in
systems having some given decay of correlations.

\begin{lemma}
\label{uno}Let $B(x_{0},r_{k})$ be a sequence of balls with decreasing
radius centered in $x_{0}$, let $A_{k}=T^{-k}(B(x_{0},r_{k}))$ and let us
write $A_{-1}=X$. If $(X,T,\mu )$ is a system satisfying definition \ref{sup}
then when $k>j>0$%
\begin{equation}
\mu (A_{k}\cap A_{j})\leq \mu (A_{k-1})\mu (A_{j-1})+\frac{\Phi (k-j)}{%
(r_{k-1}-r_{k})(r_{j-1}-r_{j})}.  \label{ball}
\end{equation}
\end{lemma}

\begin{proof}
Let $\phi _{k}$ be a Lipschitz function such that $\phi _{k}(x)=1$ for all $%
x\in B(x_{0},r_{k})$, $\phi _{k}(x)=0$ if $x\notin B(x_{0},r_{k-1})$ and $%
||\phi _{k}||\leq \frac{1}{r_{k-1}-r_{k}}$ (such functions can be easily
constructed as $\phi _{k}(x)=h(d(x_{0},x))$ where $h$ is a suitable
piecewise linear Lipshitz function $\mathbb{R\rightarrow }[0,1]$). Let $%
k>j>0 $. Since $\mu $ is preserved%
\begin{equation*}
\mu (A_{k}\cap A_{j})=\mu (T^{-k+j}(B(x_{0},r_{k}))\cap B(x_{0},r_{j}))\leq
\int \phi _{k}\circ T^{k-j}\phi _{j}d\mu \leq
\end{equation*}

\begin{equation*}
\leq \int \phi _{k}d\mu \int \phi _{j}d\mu +\left\vert \left\vert \phi
_{k}\right\vert \right\vert \left\vert \left\vert \phi _{j}\right\vert
\right\vert \Phi (k-j)\leq \mu (A_{k-1})\mu (A_{j-1})+\left\vert \left\vert
\phi _{k}\right\vert \right\vert \left\vert \left\vert \phi _{j}\right\vert
\right\vert \Phi (k-j)
\end{equation*}%
which gives the statement.
\end{proof}

In the following we will prove that if a decreasing sequence of sets $S_{k}$
are rapidly mixed like in eq. \ref{ball}, then $T^{-k}(S_{k})$ covers $X$ in
a uniform way, i.e. the $S_{k}$ form a SBC sequence. The idea of the proof
of this is somewhat similar to the proof of the strong law of large numbers
using the Paley-Zygmund inequality, and to the proof of theorem 1 in \cite%
{FMP}. The idea is to estimate $E((Z_{n})^{2})$ (this will be done by using
something similar to eq. \ref{ball}) and find an upper bound which ensures
that the distribution of the possible values of $Z_{n}$ is not too far from
the average $E(Z_{n})$. After this, choosing suitable subsequences we prove
also pointwise convergence, so that $\frac{Z_{k}}{E(Z_{k})}\rightarrow
1,~\mu -a.e.$

\begin{lemma}
\label{due}Let $S_{k}$ be a decreasing sequence of measurable sets such that 
\begin{equation*}
\underset{k\rightarrow \infty }{\lim \inf }\frac{\log (\sum_{0}^{k}\mu
(S_{k}))}{\log (k)}=z>0.
\end{equation*}%
Let $A_{k}=T^{-k}(S_{k})$ and let us suppose that the system is such that
when $k>j$%
\begin{equation}
\mu (A_{k}\cap A_{j})\leq \mu (A_{k-1})\mu (A_{j-1})+k^{c_{1}}j^{c_{2}}\Phi
(k-j)  \label{mixxx}
\end{equation}%
with $\Phi $ having superpolynomial decay and $c_{1},c_{2}\geq 0$. Then
posing $Z_{k}(x)=\sum_{0}^{k}1_{A_{i}}(x)$ we have $\frac{Z_{k}}{E(Z_{k})}%
\rightarrow 1$ in the $L^{2}$ norm and almost everywhere.
\end{lemma}

\begin{proof}
Let us estimate $E((Z_{n})^{2}).$ We have that%
\begin{equation*}
E((Z_{n})^{2})=\sum_{k=1}^{n}\mu (A_{k})+2\sum_{k,j\leq n,k>j}\mu (A_{k}\cap
A_{j}).
\end{equation*}%
Now let $0<\alpha <\frac{z}{2}$ and let us estimate the second summand on
the right side by dividing it into two parts in the following way%
\begin{equation*}
\sum_{k,j\leq n,k>j}\mu (A_{k}\cap A_{j})\leq \sum_{k,j\leq
n,k>j,k<j+n^{\alpha }}\mu (A_{k}\cap A_{j})+\sum_{k,j\leq n,k\geq
j+n^{\alpha }}\mu (A_{k}\cap A_{j}).
\end{equation*}%
The first sum can be estimated as follows: 
\begin{equation*}
\sum_{k,j\leq n,k>j,k<j+n^{\alpha }}\mu (A_{k}\cap A_{j})\leq n^{\alpha
}E(Z_{n}).
\end{equation*}%
In the second one we use equation \ref{mixxx} and we obtain

\begin{equation}
\sum_{k,j\leq n,k\geq j+n^{\alpha }}\mu (A_{k}\cap A_{j})\leq \sum_{k,j\leq
n,k\geq j+n^{\alpha }}\mu (A_{k-1})\mu (A_{j-1})+n^{c_{1}+c_{2}}\Phi
(n^{\alpha })\leq
\end{equation}%
\begin{equation*}
\leq \frac{1}{2}(E(Z_{n}))^{2}+n^{2+c_{1}+c_{2}}\Phi (n^{\alpha }).
\end{equation*}%
Hence 
\begin{equation*}
\sum_{k,j\leq n,k>j}\mu (A_{k}\cap A_{j})\leq n^{\alpha }E(Z_{n})+\frac{1}{2}%
(E(Z_{n}))^{2}+n^{2+c_{1}+c_{2}}\Phi (n^{\alpha })
\end{equation*}%
and%
\begin{equation*}
E((Z_{n})^{2})\leq (2n^{\alpha
}+1)E(Z_{n})+(E(Z_{n}))^{2}+2n^{2+c_{1}+c_{2}}\Phi (n^{\alpha }).
\end{equation*}%
Now, let us remark that 
\begin{equation*}
E((Z_{n}-E(Z_{n}))^{2})=E((Z_{n})^{2})-(E(Z_{n}))^{2}
\end{equation*}%
hence 
\begin{equation}
E((Z_{n}-E(Z_{n}))^{2})\leq (2n^{\alpha }+1)E(Z_{n})+2n^{2+c_{1}+c_{2}}\Phi
(n^{\alpha }).
\end{equation}%
Now we want to compare $Z_{n}$ with its average $E(Z_{n}).$ For this we
consider 
\begin{equation*}
Y_{n}=\frac{Z_{n}}{E(Z_{n})}-1=\frac{Z{_{n}}-{E(Z_{n})}}{E(Z_{n})}.
\end{equation*}%
When $Y_{n}=0$, $Z_{n}=E(Z_{n})$. By the above results 
\begin{equation}
E((Y_{n})^{2})\leq \frac{(2n^{\alpha }+1)E(Z_{n})+2n^{2+c_{1}+c_{2}}\Phi
(n^{\alpha })}{(E(Z_{n}))^{2}},
\end{equation}%
since $\alpha <\frac{z}{2}$ and $2n^{2+c_{1}+c_{2}}\Phi (n^{\alpha
})\rightarrow 0,$ then $\underset{n\rightarrow \infty }{lim}E((Y_{n})^{2})=0$%
. This proves that $\frac{Z_{n}}{E(Z_{n})}\rightarrow 1$ in $L^{2}$.

By this it is easy to see that there is a subsequence of $Y_{n}$ that
converges a.e. to $0$, but we want to prove that the whole sequence
converges to $0$ a.e. We begin by considering%
\begin{equation}
n_{k}=\inf \{{n:E(Z_{n})\geq k^{2}\}}
\end{equation}%
and show that $Y_{n_{k}}\rightarrow 0$ a.e. We remark that since $\Phi $ has
a superpolynomial decay 
\begin{equation}
\sum_{n}n^{2+c_{1}+c_{2}}\frac{\Phi (n^{\alpha })}{(E(Z_{n}))^{2}}<\infty .
\label{superpol}
\end{equation}%
Consider a small $\epsilon >0$, by definition of $z$ and by the fact that $%
\mu (A_{i})<1$ we have that if $k$ is big enough $(k+1)^{2}\geq
E(Z_{n_{k}})\geq (n_{k})^{z-\epsilon }.$ Hence $n_{k}\leq (k+1)^{\frac{2}{%
z-\epsilon }}\leq (2k)^{\frac{2}{z-\epsilon }}$ and so $\frac{%
(2n_{k}^{\alpha }+1)E(Z_{n_{k}})}{(E(Z_{n_{k}}))^{2}}=\frac{2n_{k}^{\alpha
}+1}{E(Z_{n_{k}})}\leq \frac{2(2k)^{\frac{2\alpha }{z-\epsilon }}+1}{k^{2}}$%
. Since $\alpha <\frac{z}{2},$ we can suppose $\epsilon $ to be so small
that $2\alpha <z-\epsilon $ and hence%
\begin{equation}
\frac{2\alpha }{z-\epsilon }-2<-1.  \label{zetamezzi}
\end{equation}%
This implies that $\sum E((Y_{n_{k}})^{2})<\infty $, which in turn implies
that $Y_{n_{k}}\rightarrow 0$ a.e. and $\frac{Z_{n_{k}}}{E(Z_{n_{k}})}%
\rightarrow 1$ almost everywhere.

Now, if $n_{k}\leq n\leq n_{k+1}$

\begin{equation*}
\frac{Z_{n}}{E(Z_{n})}\leq \frac{Z_{n_{k+1}}}{E(Z_{n_{k}})}=\frac{Z_{n_{k+1}}%
}{E(Z_{n_{k+1}})}\frac{E(Z_{n_{k+1}})}{E(Z_{n_{k}})}\leq \frac{Z_{n_{k+1}}}{%
E(Z_{n_{k+1}})}\frac{(k+2)^{2}}{k^{2}}
\end{equation*}

and

\begin{equation*}
\frac{Z_{n}}{E(Z_{n})}\geq \frac{Z_{n_{k}}}{E(Z_{n_{k+1}})}=\frac{Z_{n_{k}}}{%
E(Z_{n_{k}})}\frac{E(Z_{n_{k}})}{E(Z_{n_{k+1}})}\geq \frac{Z_{n_{k}}}{%
E(Z_{n_{k}})}\frac{k^{2}}{(k+2)^{2}}.
\end{equation*}

then we have $\underset{n\rightarrow \infty }{\lim }\frac{Z_{n}}{E(Z_{n})}%
=1, $ $\mu -$almost everywhere.
\end{proof}

\begin{remark}
\label{remarkbc}We remark that in the above proof, the key point to ensure
that $\sum E((Y_{n_{k}})^{2})<\infty $ (and then have a.e. convergence) are
equations \ref{superpol} and \ref{zetamezzi} this implies that if $\Phi $
has not superpolinomial decay, but a polinomial decay fast enough that $%
\sum_{n}n^{2+c_{1}+c_{2}}\frac{\Phi (n^{\alpha })}{(E(Z_{n}))^{2}}<\infty $
for some $\alpha <\frac{z}{2}$ then $\underset{n\rightarrow \infty }{lim}%
\frac{Z_{n}}{E(Z_{n})}=1,$ $\mu -$almost everywhere. Then the lemma also
holds for some system such that $\Phi $ has (rapid enough) polinomial decay.
\end{remark}

Now we use the above result to conclude the equality between the hitting
time indicator and dimension. The second part of the proof is similar to the
proof of Theorem 2.4 in \cite{GK}.

\begin{proof}
\emph{(of Thm. \ref{maine})} Let us prove $\overline{R}(x,x_{0})\leq d_{\mu
}(x_{0})$ for almost each $x.$ We recall that this implies \underline{$R$}$%
(x,x_{0})\leq d_{\mu }(x_{0})$ and the opposite inequalities come from
theorem \ref{GAN}. Let us consider $0<\beta <\frac{1}{d_{\mu }(x_{0})}$ and
a sequence $r_{k}=k^{-\beta }$ (we remark that if the result is proved for
such a subsequence, hence it holds for all subsequences, see lemma 4.2 in 
\cite{GKP}). Then for each small $\epsilon <\beta ^{-1}-d_{\mu }(x_{0})$,
eventually $\mu (B(x_{0},r_{k}))\geq (r_{k})^{d_{\mu }(x_{0})+\epsilon
}=k^{-\beta (d_{\mu }(x_{0})+\epsilon )}$ and if $k$ is big enough $%
\sum_{0}^{k}\mu (B(x_{0},r_{i}))\geq Ck^{1-\beta (d_{\mu }(x_{0})+\epsilon
)},$ since $\epsilon $ is arbitrary we have 
\begin{equation*}
\underset{k\rightarrow \infty }{\lim \inf }\frac{\log (\sum_{0}^{k}\mu
(B(x_{0},r_{i})))}{\log (k)}\geq 1-\beta d_{\mu }(x_{0})>0.
\end{equation*}%
Moreover, $r_{k-1}-r_{k}\sim k^{-\beta -1}.$ Hence we can apply Lemma \ref%
{uno} and \ref{due} to the sequence $B(x_{0},r_{k})$ and obtain that for
such a sequence $\underset{n\rightarrow \infty }{lim}\frac{Z_{n}}{E(Z_{n})}%
=1,$ $\mu -$almost everywhere.

Let us consider $\epsilon ^{\prime }>0$ and $\beta $ as above, near to $%
\frac{1}{d_{\mu }(x_{0})}$, such that $\beta (d_{\mu }(x_{0})+\epsilon
^{\prime })>1$. Moreover let us consider $\varepsilon >0$ so small that $%
\beta (d_{\mu }(x_{0})+\varepsilon )<1$ and $\beta (d_{\mu }(x_{0})+\epsilon
^{\prime })-\frac{1-\beta (d_{\mu }(x_{0})-\varepsilon )}{1-\beta (d_{\mu
}(x_{0})+\varepsilon )}>0$. Let us consider $x$ such that $\overline{R}%
(x,x_{0})>d_{\mu }(x_{0})+\epsilon ^{\prime }$, then for infinitely many $n$%
, $\tau _{n^{-\beta }}(x,x_{0})>n^{\beta (d_{\mu }(x_{0})+\epsilon ^{\prime
})}$. Then 
\begin{equation*}
x\notin \cup _{0\leq i\leq n^{\beta (d_{\mu }(x_{0})+\epsilon ^{\prime
})}}T^{-i}(B(x_{0},n^{-\beta }))
\end{equation*}%
and in particular 
\begin{equation*}
x\notin \cup _{n\leq i\leq n^{\beta (d_{\mu }(x_{0})+\epsilon ^{\prime
})}}T^{-i}(B(x_{0},n^{-\beta }))\supset \cup _{n\leq i\leq n^{\beta (d_{\mu
}(x_{0})+\epsilon ^{\prime })}}T^{-i}(B(x_{0},i^{-\beta })),
\end{equation*}%
which implies that there is a sequence $n_{i}$ such that $%
Z_{n_{i}}(x)=Z_{n_{i}^{\beta (d_{\mu }(x_{0})+\epsilon ^{\prime })}}(x)$ for
each $i.$ Now let us consider $E(Z_{n_{i}})$ and $E(Z_{n_{i}^{\beta (d_{\mu
}(x_{0})+\epsilon ^{\prime })}}).$ Since the local dimension at $x_{0}$ is $%
d_{\mu }(x_{0})$,\ when $i$ is big enough 
\begin{equation*}
i^{-\beta (d_{\mu }(x_{0})+\varepsilon )}<\mu (B(x_{0},i^{-\beta
}))<i^{-\beta (d_{\mu }(x_{0})-\varepsilon )}
\end{equation*}%
then there are constants $k_{1}$ and $k_{2}$ such that when $n$ is big
enough $k_{1}n^{1-\beta (d_{\mu }(x_{0})+\varepsilon
)}<E(Z_{n})<k_{2}n^{1-\beta (d_{\mu }(x_{0})-\varepsilon )}$. From this we
have that if $i$ is big enough

\begin{eqnarray*}
\frac{E(Z_{n_{i}})}{E(Z_{n_{i}^{\beta (d_{\mu }(x_{0})+\epsilon ^{\prime
})}})} &\leq &\frac{k_{2}n_{i}^{1-\beta (d_{\mu }(x_{0})-\varepsilon )}}{%
k_{1}n_{i}^{\beta (d_{\mu }(x_{0})+\epsilon ^{\prime })(1-\beta (d_{\mu
}(x_{0})+\varepsilon ))}}= \\
&=&\frac{k_{2}}{k_{1}}n_{i}^{(1-\beta (d_{\mu }(x_{0})-\varepsilon ))-\beta
(d_{\mu }(x_{0})+\epsilon ^{\prime })(1-\beta (d_{\mu }(x_{0})+\varepsilon
))}.
\end{eqnarray*}%
By the assumptions on $\varepsilon ,$ $(1-\beta (d_{\mu }(x_{0})-\varepsilon
))-\beta (d_{\mu }(x_{0})+\epsilon ^{\prime })(1-\beta (d_{\mu
}(x_{0})+\varepsilon ))=(1-\beta (d_{\mu }(x_{0})+\varepsilon ))(\frac{%
1-\beta (d_{\mu }(x_{0})-\varepsilon )}{1-\beta (d_{\mu }(x_{0})+\varepsilon
)}-\beta (d_{\mu }(x_{0})+\epsilon ^{\prime }))<0,$ hence%
\begin{equation*}
\underset{i\rightarrow \infty }{\lim }\frac{E(Z_{n_{i}})}{E(Z_{n_{i}^{\beta
(d_{\mu }(x_{0})+\epsilon ^{\prime })}})}=0.
\end{equation*}%
Since $n_{i}$ was chosen such that $Z_{n_{i}}(x)=Z_{n_{i}^{\beta (\overline{d%
}_{\mu }(x_{0})+\epsilon ^{\prime })}}(x)$ this implies that $\frac{%
Z_{n_{i}}(x)}{E(Z_{n_{i}})}\frac{E(Z_{n_{i}^{\beta (d_{\mu }(x_{0})+\epsilon
^{\prime })}})}{Z_{n_{i}^{\beta (d_{\mu }(x_{0})+\epsilon ^{\prime })}}(x)}=%
\frac{E(Z_{n_{i}^{\beta (d_{\mu }(x_{0})+\epsilon ^{\prime })}})}{%
E(Z_{n_{i}})}\rightarrow \infty $ as $i$ increases. Then is not possible
that $\underset{n\rightarrow \infty }{\lim }\frac{Z_{n}(x)}{E(Z_{n}(x))}=1.$
This, implies that $\overline{R}(x,x_{0})>d_{\mu }(x_{0})+\epsilon ^{\prime
} $ on a zero measure set. Finally, since $\epsilon ^{\prime }$ can be
chosen to be arbitrarily small we have the statement.
\end{proof}

We end by remarking that the above lemmas \ref{uno}, \ref{due} and remark %
\ref{remarkbc} give the following consequence, which in the author opinion
is interesting by itself (by lemma \ref{due}, the proof is obtained as in
the first part of the proof of thm. \ref{maine}):

\begin{corollary}
\label{XXX2} If the system has decay of correlation given by the function $%
\Phi $ (see definition \ref{sup}) and $B(x_{0},r_{i})$ is a sequence of
decreasing balls such that

i) $\underset{k\rightarrow \infty }{\lim \inf }\frac{\log (\sum_{0}^{k}\mu
(B(x_{0},r_{i})))}{\log (k)}=z>0$,

ii) $\underset{n\rightarrow \infty }{\lim \inf }$ $\frac{\log (r_{n-1}-r_{n})%
}{\log (n)}=c\in \mathbb{R}$

iii) $\sum_{n}n^{2-2c+\epsilon }\frac{\Phi (n^{\alpha })}{(\sum_{0}^{k}\mu
(B(x_{0},r_{i})))^{2}}<\infty $ for some $\alpha <\frac{z}{2}$ and $\epsilon
>0$

then $B(x_{0},r_{i})$ has the SBC property.
\end{corollary}

\end{document}